\documentclass{amsart} 
\usepackage{graphicx}
\usepackage{amsmath}
\usepackage{amssymb}
\usepackage{amsthm}
\newtheorem{thm}{Theorem}

\newtheorem{lem}{Lemma} 
\newtheorem{definition}{Definition}
\newtheorem{rem}{Remark}
\newtheorem{example}{Example}

\begin{document}

\title[Boundary Continuity]{Boundary Continuity of Nonparametric Prescribed Mean Curvature Surfaces}
\date{\today}
\author{Mozhgan ``Nora'' Entekhabi }
\address{Department of Mathematics \\ Florida A \& M University \\ Tallahassee, FL 32307}
\email{mozhgan.entekhabi@famu.edu}
\author{Kirk E. Lancaster}
\address{Wichita, Kansas 67226}
\thanks{The second author wishes to thank the Institute of Mathematics at Academia Sinica for support during part of this investigation and to 
thank Professor Fei-tsen Liang for helpful discussions.}
\email{redwoodsrunner@gmail.com}

\subjclass[2010]{Primary: 35J67; Secondary: 35J93, 53A10}
\keywords{prescribed mean curvature, Dirichlet problem, boundary continuity}
 
\def\Real{{\rm I\hspace{-0.2em}R}}
\def\Natural{{\rm I\hspace{-0.2em}N}}
\newcommand\myeq{\mathrel{\overset{\makebox[0pt]{\mbox{\normalfont\tiny\sffamily def}}}{=}}}

\begin{abstract}
We investigate the boundary behavior of variational solutions of Dirichlet problems for prescribed mean curvature equations at smooth boundary points where 
certain boundary curvature conditions are satisfied (which preclude the existence of local barrier functions).     
We prove that if the Dirichlet boundary data $\phi$  is continuous at such a point (and possibly nowhere else), then the solution of the variational 
problem is continuous at this point.
\end{abstract}

\maketitle

\section{Introduction}
Let $\Omega$  be a locally Lipschitz domain in $\Real^{2}$  and define $Nf = \nabla \cdot Tf = {\rm div}\left(Tf\right),$  where  $f\in C^{2}(\Omega)$  and 
$Tf= \frac{\nabla f}{\sqrt{1+\left|\nabla f\right|^{2}}}.$   
Let $H\in C^{1}(\overline{\Omega})$  satisfy the condition 
\begin{equation}
\label{(16.60)}
\left|\int_{\Omega} H\eta \ dx\right| \le \frac{1}{2} \int_{\Omega} |D\eta|\ dx \ \ \ \ {\rm for \ all \ } \eta\in C^{1}_{0}(\Omega)
\end{equation}
(e.g. \cite[(16.60)]{GT}). 
Here and throughout the paper, we adopt the sign convention that the curvature of $\Omega$  is nonnegative when $\Omega$  is convex.  
Consider the Dirichlet problem
\begin{eqnarray}
\label{eq:D}
Nf & = & 2H  \mbox{  \ in \ } \Omega    \\
f & = & \phi  \mbox{ \ on \ } \partial \Omega.
 \label{bc:D}
\end{eqnarray}
We wish to understand the boundary behavior of a solution of (\ref{eq:D})-(\ref{bc:D}).  

If $\partial\Omega$  is smooth and the curvature of $\partial\Omega$  is greater than or equal to $2|H|$  at each point of $\partial\Omega,$  
then this problem is well-posed. 
If $\phi\in C^{0}(\partial\Omega),$  then there exists a unique solution of (\ref{eq:D})-(\ref{bc:D}) in $C^{2}(\Omega)\cap C^{0}(\overline{\Omega})$  
and $f=\phi$  on $\partial\Omega$  (\cite[Theorem 16.10]{GT}).  
On the other hand, one can choose a distinguished point ${\mathcal O}\in\partial\Omega$  and use the ``gliding hump'' construction as in \cite{Lan:89} 
and \cite[Theorem 3]{LS1}, in conjunction with \cite[Theorem 2]{NoraKirk1},   
to prove that there exist $\phi\in L^{\infty}(\partial\Omega)$  such that the (unique) variational solution $f$  of (\ref{eq:D})-(\ref{bc:D}) is 
in $C^{2}(\Omega)$  but $f$  is discontinuous at ${\mathcal O}$  and none of the radial limits of $f$  at ${\mathcal O}$  exist.  

When the appropriate geometric conditions (e.g. convexity of the domain in the case of the minimal surface equation on $\Real^{2}$) are not satisfied, then 
the Dirichlet problem is ill-posed and a classical solution of (\ref{eq:D})-(\ref{bc:D}) may not exist (see, for example, \cite{JS1968} and \cite[\S 406]{NitscheBook}).   
Interest in determining sufficient conditions for the existence of classical solutions of (\ref{eq:D})-(\ref{bc:D}) is long standing and 
one method is to impose ``smallness'' conditions on the Dirichlet boundary data $\phi.$
When $H\equiv 0$  in (\ref{eq:D}), an early result is A. Korn's classic 1909 paper (\cite{Korn}); J.C.C. Nitsche discusses some of the history 
of this problem for the minimal surface equation in \cite[\S 285 \& \S 412]{NitscheBook}.   
G. Williams (\cite{Williams1984,Williams1985}), C.P. Lau (\cite{Lau}), K. Hayasida and M. Nakatani (\cite{HayasidaNakatani}), 
M. Bergner (\cite{Bergner}), J. Ripoll and F. Tomi (\cite{RT}) and 
many others have investigated limiting $\phi$  in order to prove a classical solution exists. 

Rather than imposing a ``smallness'' condition on the Dirichlet data $\phi$  and trying to obtain a classical solution of (\ref{eq:D})-(\ref{bc:D}),  
we wish to impose a (local) condition on the curvature of the domain, place no restrictions on the Dirichlet data $\phi\in L^{\infty}(\partial\Omega)$  and   
prove that the variational solution $f$  extends to be continuous at a point $\mathcal{O}$  of $\partial\Omega$  (or on an open subset $\Gamma$  of $\partial\Omega$)  
in the sense that $f\in C^{0}\left(\Omega\cup \{\mathcal{O}\}\right)$  (or $f\in C^{0}\left(\Omega\cup \Gamma\right)$) 
when $\phi$  is continuous at $\mathcal{O}$  (or on $\Gamma$).  
In general, no classical solution may exist and the variational solution is the best approximation to a classical solution.

We shall assume that $\partial\Omega$  is smooth, ${\mathcal O}\in\partial\Omega$  is a distinguished point, the curvature $\Lambda$  of $\partial\Omega$  satisfies 
\begin{equation}
\label{Condition1}
\Lambda({\bf x}) < -2|H({\bf x})| \ \ \ \ {\rm for} \ {\bf x}\in\partial\Omega, \ |{\bf x}-{\mathcal O}|<\delta
\end{equation}
for some $\delta>0$  and  $\phi\in L^{\infty}(\partial\Omega).$  
In Theorem \ref{Thm1}, we prove that the (unique) variational solution $f\in C^{2}(\Omega)$  of (\ref{eq:D})-(\ref{bc:D}) is 
continuous at ${\mathcal O}$  if $\phi$  is continuous at ${\mathcal O}.$  
In Theorem \ref{Thm2}, we assume $H\equiv 0,$  relax (\ref{Condition1}) slightly and obtain the same result as in Theorem \ref{Thm1}.
In Theorem \ref{Thm3}, we prove that the radial limits $Rf(\cdot)$  of $f$  exist at ${\mathcal O}$   if $\phi$  is discontinuous at ${\mathcal O}$ 
(even if $\phi$  does not have one-sided limits at ${\mathcal O}$).  

The idea of the proof is to describe the graph of $f$  parametrically in isothermal coordinates, prove that it is uniformly continuous on its (open) 
parameter domain and therefore extends uniquely to a continuous function on the closure of the parameter domain.  
In the case of Theorems \ref{Thm1} and \ref{Thm2}, this is equivalent to proving that $f$  is uniformly continuous on the intersection of $\Omega$  and 
an open neighborhood $U$  of ${\mathcal O}$  and therefore $f$  extends uniquely to a continuous function on $\overline{\Omega\cap U}.$
Of course, just as in \cite[Theorem 4.2]{Bour}, \cite{Lin} and \cite{Simon},  this extension of $f$  needs not equal $\phi.$  
Differences between our results and those, for example, in the papers by Bourni, Lin and Simon are the additional requirements imposed on the domain or the boundary 
data; \cite{Bour}  requires the graph of $\phi$  to be a specific type of limit of the graphs of $C^{1,\alpha}$  functions, 
\cite{Lin}  requires $\phi$  to be Lipschitz and \cite{Simon} requires $H\equiv 0,$  $\partial\Omega$  to be $C^{4}$  and  $\phi$  to be Lipschitz. 
In Example \ref{Example1}, we present an illustration of the value of symmetry and prove that $f$  can be continuous at ${\mathcal O}$  even if $\phi$  
is not; this illustration uses the same domain as that mentioned in \cite{Simon}.     

In the remaining case where $-2|H({\mathcal O})|\le \Lambda({\mathcal O})<2|H({\mathcal O})|,$  the behavior of the variational solution at 
${\mathcal O}$  is unknown.  The Dirichlet problem (\ref{eq:D})-(\ref{bc:D}) 
does not have a classical solution in $C^{2}(\Omega)\cap C^{0}(\overline{\Omega})$  for all $\phi\in C^{\infty}(\partial\Omega)$  (\cite[Corollary 14.13]{GT}) 
and the standard gliding hump argument does not work in this situation.  
On the other hand, not all of the comparison functions needed here are available and so the conclusions of Theorems \ref{Thm1} and \ref{Thm2} may not hold. 
Moonies (\cite{FinnMoon,Liang92}, see also \cite{Liang96}) are not bounded below but do illustrate some properties of this case, 
with $H\equiv 1,$  $\Lambda\equiv -2$  on one component of the boundary and $\Lambda\equiv \frac{1}{R},$  $\frac{1}{2}<R<1,$  on the other component. 
As in \cite{JL}, these cases illustrate the strong differences between uniformly elliptic (genre zero) and prescribed mean curvature (genre two) equations 
(see also Remark \ref{Remark2}); 
for Laplace's equation in $\Real^{2},$  for example, boundary curvature would play no role in the solvability of Dirichlet problems and 
the gliding hump construction could always be used, exactly as in the first case above. 
\vspace{3mm}

\begin{thm}
\label{Thm1}
Suppose $\Omega$  is a locally Lipschitz  domain in $\Real^{2},$  $\Gamma$  is a $C^{2,\lambda}$  open subset of $\partial\Omega$  for some $\lambda\in (0,1),$    
the  curvature $\Lambda({\bf x})$   of $\Gamma$  at ${\bf x}$  is less than $-2|H({\bf x})|$  for ${\bf x}\in\Gamma.$  
Suppose $\phi\in L^{\infty}(\partial\Omega),$  ${\bf y}\in \Gamma,$  
either $f$  is symmetric with respect to a line through ${\bf y}$  or $\phi$  is continuous at ${\bf y},$  and $f\in BV(\Omega)$  minimizes 
\begin{equation}
\label{AA}
J(u)=\int_{\Omega} \sqrt{1+|Du|^{2}}\ d{\bf x} + \int_{\Omega}  2Hu \ d{\bf x} + \int_{\partial\Omega} |u-\phi|\ ds
\end{equation}
for $u\in BV(\Omega).$  Then $f\in C^{0}(\Omega\cup \{ {\bf y}\}).$  
If $\phi\in C^{0}(\Gamma),$  then $f\in C^{0}(\Omega\cup\Gamma).$  
\end{thm}
\vspace{3mm}

\noindent When $H\equiv 0,$  the strict inequality (\ref{Condition1}) can be relaxed.   

\begin{thm}
\label{Thm2}
Suppose $\Omega$  is a locally Lipschitz  domain in $\Real^{2},$  $\Gamma$  is a $C^{2,\lambda}$  open subset of $\partial\Omega$  
for some $\lambda\in (0,1),$  and the  curvature $\Lambda$  of $\Gamma$  is nonpositive and vanishes, at most, at a finite number of points of $\Gamma.$  
For each point ${\bf x}_{0}\in \Gamma$  at which $\Lambda({\bf x}_{0})=0,$  suppose there exist $C>0$  and $\delta>0$  such that 
\begin{equation}
\label{c-bound}
|\Lambda({\bf x})|\ge C|{\bf x}-{\bf x}_{0}|^{\lambda} \ \ {\rm for} \ \  {\bf x}\in \Gamma, \ |{\bf x}-{\bf x}_{0}|<\delta.
\end{equation} 
Suppose $\phi\in L^{\infty}(\partial\Omega),$  ${\bf y}\in \Gamma,$  
either $f$  is symmetric with respect to a line through ${\bf y}$  or $\phi$  is continuous at ${\bf y},$  and $f\in BV(\Omega)$  minimizes 
\begin{equation}
\label{JJ}
J(u)=\int_{\Omega} \sqrt{1+|Du|^{2}}\ d{\bf x} + \int_{\partial\Omega} |u-\phi| \ ds
\end{equation}
for $u\in BV(\Omega).$  Then $f\in C^{0}(\Omega\cup \{ {\bf y}\}).$  
If $\phi\in C^{0}(\Gamma),$  then $f\in C^{0}(\Omega\cup\Gamma).$  
\end{thm} 
\vspace{3mm}

\begin{figure}[htb]
\centerline{
\includegraphics[width=1.5in]{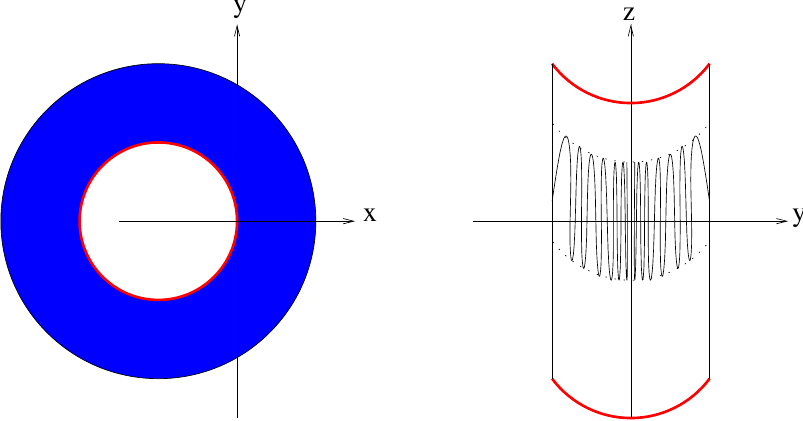}}
\caption{$\Omega$  and part of the graph of $\phi$  \label{Cat}}
\end{figure}

\begin{example}
\label{Example1}
Let $\Omega=\{(x,y)\in \Real^{2} : 1< (x+1)^{2}+y^{2}<\cosh^{2}(1)\}$  and $\phi(x,y)=\sin\left(\frac{\pi}{x^2+y^2}\right)$  for $(x,y)\neq (0,0)$  
(see Figure \ref{Cat} for a rough illustration of the graph of $\phi$).  Set ${\mathcal O}=(0,0)$  and   $H\equiv 0.$ 
Let $f\in C^{2}(\Omega)$  minimize (\ref{JJ}) over $BV(\Omega).$   
Then Theorem \ref{Thm1} (with ${\bf y}={\mathcal O}$) implies $f\in C^{0}\left(\overline{\Omega}\right),$  even though $\phi$  has no limit at ${\mathcal O}.$  
\end{example}

\noindent We shall prove the following theorem on the existence and behavior of the radial limits of variational solutions of (\ref{eq:D})-(\ref{bc:D}) 
and use this to prove Theorems \ref{Thm1} and \ref{Thm2}.  
At a point ${\bf y}\in\partial\Omega,$  we let $\alpha({\bf y})$  and $\beta({\bf y})$  be the angles which the tangent rays to $\partial\Omega$  at 
${\bf y}$  make with the positive $x-$axis such that 
\[
\{ {\bf y}+r(\cos\theta,\sin\theta) : 0<r<\epsilon(\theta), \alpha({\bf y})<\theta<\beta({\bf y}) \} \subset \Omega\cap B_{\delta}({\bf y}) 
\]
for some $\delta>0$  and some function $\epsilon(\cdot):(\alpha({\bf y}),\beta({\bf y}))\to (0,\delta).$ 
If $\partial\Omega$  is smooth at ${\bf y},$  $\beta({\bf y})=\alpha({\bf y})+\pi.$  
For $\theta\in (\alpha({\bf y}),\beta({\bf y})),$  $Rf(\theta,{\bf y})=\lim_{r\downarrow 0}f({\bf y}+r(\cos\theta,\sin\theta))$  if this limit exists. 
Also $Rf(\alpha({\bf y}),{\bf y})=\lim_{\Gamma_{1}({\bf y})\ni {\bf x}\to {\bf y}} f({\bf x})$  if this limit exists and 
$Rf(\beta({\bf y}),{\bf y})=\lim_{\Gamma_{2}({\bf y})\ni {\bf x}\to {\bf y}} f({\bf x})$  if this limit exists, 
where $\partial \Omega \cap B_{\delta}({\bf y})\setminus \{{\bf y}\}$ consists of disjoint, open arcs  $\Gamma_{1}({\bf y})$  and $\Gamma_{2}({\bf y})$  
whose tangent rays approach the rays $\theta= \alpha({\bf y})$  and $\theta= \beta({\bf y})$ respectively, as the point ${\bf y}$ is approached. 
\vspace{5mm}

\begin{thm}
\label{Thm3}
Suppose $\Omega$  is a locally Lipschitz  domain in $\Real^{2},$  $\Gamma$  is a $C^{2,\lambda}$  open subset of $\partial\Omega$  
for some $\lambda\in (0,1).$  Suppose either (a) $H\equiv 0$  in a neighborhood of $\Gamma,$   
the curvature $\Lambda$  of $\Gamma$  is nonpositive and vanishes at only a finite number of points of $\Gamma,$  at each of which (\ref{c-bound}) holds,   
or  (b) the  curvature $\Lambda({\bf x})$   of $\Gamma$  at ${\bf x}$  is less than $-2|H({\bf x})|$  for ${\bf x}\in\Gamma.$  
Let $f\in C^{2}(\Omega)\cap L^{\infty}(\Omega)$  satisfy $Nf=2H$  in $\Omega.$  
Suppose that ${\bf y}\in\partial\Omega.$     Then the limits 
\[
\lim_{\Gamma_{1}({\bf y})\ni {\bf x}\to {\bf y}} f({\bf x})=z_{1} \ {\rm and} \ 
\lim_{\Gamma_{2}({\bf y})\ni {\bf x}\to {\bf y}} f({\bf x})=z_{2}
\]
exist, $Rf(\theta,{\bf y})$  exists for each $\theta\in [\alpha({\bf y}),\beta({\bf y})],$
$Rf(\cdot,{\bf y})\in C^{0}([\alpha({\bf y}),\beta({\bf y})]),$  and $Rf(\cdot,{\bf y})$ behaves in one of the following ways:

\noindent (i) $Rf(\cdot,{\bf y})=z_{1}$  is a constant function  and $f$  is continuous at ${\bf y}.$  

\noindent (ii) There exist $\alpha_{1}$ and $\alpha_{2}$ so that $\alpha({\bf y}) \leq \alpha_{1}
< \alpha_{2} \leq \beta({\bf y}),$  $Rf=z_{1}$  on $[\alpha({\bf y}), \alpha_{1}],$  $Rf=z_{2}$  on $[\alpha_{2}, \beta({\bf y})]$
and $Rf$ is strictly increasing (if $z_{1}<z_{2}$) or strictly decreasing  (if $z_{1}>z_{2}$)  on $[\alpha_{1}, \alpha_{2}].$    
\end{thm}

\section{Proofs}

Let  $Q$  be the operator on $C^{2}(\Omega)$  given by 
\begin{equation}
\label{Q}
Qf({\bf x}) \myeq Nf({\bf x}) - 2H({\bf x}), \ \ \ \ {\bf x}\in\Omega. 
\end{equation}
Let $\nu$  be the exterior unit normal to $\partial\Omega,$  defined almost everywhere on $\partial\Omega.$    
We assume that for almost every ${\bf y}\in\partial\Omega,$  there is a continuous extension $\hat\nu$  of $\nu$  to a neighborhood of ${\bf y}.$ 

\begin{definition}
\label{Bernstein+} 
Given a locally Lipschitz domain $\Omega,$ a \underline{\bf upper Bernstein pair} $\left(U^{+},\psi^{+}\right)$  for a curve $\Gamma\subset \partial\Omega$  
and a function $H$  in (\ref{Q}) is a $C^{1}$  domain $U^{+}$  and a function $\psi^{+}\in C^{2}(U^{+})\cap C^{0}\left(\overline{U^{+}}\right)$  such that 
$\Gamma\subset\partial U^{+},$  $\nu$  is the exterior unit normal to $\partial U^{+}$  at each point of $\Gamma$  
(i.e. $U^{+}$  and $\Omega$  lie on the same side of $\Gamma$), $Q\psi^{+}\le 0$  in $U^{+},$  and $T\psi^{+}\cdot\nu=1$  almost everywhere on an open subset 
of $\partial U^{+}$  containing $\overline{\Gamma}$  in the same sense as in  \cite{CF:74a}; that is, for almost every ${\bf y}\in\Gamma,$ 
\begin{equation}
\label{zero}
\lim_{U^{+}\ni {\bf x}\to {\bf y}} \frac{\nabla \psi^{+}({\bf x})\cdot \hat\nu({\bf x})}{\sqrt{1+|\nabla \psi^{+}({\bf x})|^{2}}} = 1.
\end{equation} 
\end{definition}

\begin{definition}
\label{Bernstein-} 
Given a domain $\Omega$  as above, a \underline {\bf lower Bernstein pair} $\left(U^{-},\psi^{-}\right)$  for a curve $\Gamma\subset \partial\Omega$  
and a function $H$  in (\ref{Q}) is a $C^{1}$  domain $U^{-}$  and a function $\psi^{-}\in C^{2}(U^{-})\cap C^{0}\left(\overline{U^{-}}\right)$  
such that $\Gamma\subset\partial U^{-},$  $\nu$  is the exterior unit normal to $\partial U^{-}$  at each point of $\Gamma$  
(i.e. $U^{-}$  and $\Omega$  lie on the same side of $\Gamma$), 
$Q\psi^{-}\ge 0$  in $U^{-},$  and $T\psi^{-}\cdot\nu=-1$  almost everywhere on an open subset of $\partial U^{-}$  containing $\overline{\Gamma}$  
in the same sense as in  \cite{CF:74a}.
\end{definition}

The existence of Bernstein pairs is established in Lemmas \ref{Prop2} and \ref{Prop3}.

\begin{lem}
\label{Prop1}
Let $a<b,$  $\lambda\in (0,1),$  $\psi\in C^{2,\lambda}([a,b])$  and  
$
\Gamma=\{(x,\psi(x))\in \Real^{2} : x\in [a,b]\}
$
such that 
$\psi'(x)< 0$  for $x\in [a,b],$  $\psi''(x)<0$  for $x\in [a,b]\setminus J,$    there exist $C_{1}>0$  and $\epsilon_{1}>0$  such that 
if $\bar x\in J$  and $|x-\bar x|<\epsilon_{1},$  then $\psi''(x)\le -C_{1}|x-\bar x|^{\lambda},$  where $J$  is a finite subset of $(a,b).$
Then there exists an open set $U\subset\Real^{2}$  with  $\Gamma\subset \partial U$  and a function 
$h\in C^{2}\left(U\right)\cap C^{0}\left(\overline{U}\right)$  such that  $\partial U$  is a closed, $C^{2,\lambda}$  curve, 
$\Gamma$  lies below $U$  in $\Real^{2}$  (i.e. the exterior unit normal $\nu=(\nu_{1}(x),\nu_{2}(x))$  to $\partial U$ satisfies $\nu_{2}(x)<0$  
for $a\le x\le b$),   $Nh=0$  in $U$  and $Th\cdot\nu=1$  almost everywhere on an open subset of $\partial U^{+}$  containing $\overline{\Gamma}$  
(i.e. (\ref{zero}) holds).
\end{lem}
\vspace{3mm} 

\noindent {\bf Proof:} We may assume that $a,b>0.$  
There exists $c>b$  and $k\in C^{2,\lambda}([-c,c])$  with $k(-x)=k(x)$  for $x\in [0,c]$  such that $k(x)=-\psi(x)$  for $x\in [a,b],$  
$k''(x)>0$  for $x\in [-c,c]\setminus J,$  where $J$  is a finite set, $k''(0)>0,$  and the set 
\[
K=\{(x,k(x))\in \Real^{2} : x\in [-c,c]\}  
\]
is strictly concave (i.e. $tk(x_{1})+(1-t)k(x_{2})>k\left(tx_{1}+(1-t)x_{2}\right)$  for each $t\in (0,1)$  and 
$x_{1},x_{2}\in [-c,c]$  with $x_{1}\neq x_{2}$).  
From \cite[pp.1063-5]{EL1986B}, we can construct a domain $\Omega(K,l)$ such that $K\subset\partial\Omega(K,l)$  and 
$\Omega(K,l)$  lies below $K$  (i.e. the outward unit normal to  $\Omega(K,l)$  at $(x,k(x))$  is 
$\nu(x)=\frac{\left(-k'(x),1\right)}{\sqrt{1+\left(k'(x)\right)^{2}}};$  see \cite[Figure 4]{EL1986B}) 
and a function $F^{+}\in C^{2}\left(\Omega(K,l)\right)\cap C^{0}\left(\overline{\Omega(K,l)}\right)$  such that the function 
\[
\mu({\bf x})\myeq \frac{\left(\nabla F^{+}({\bf x}),-1\right)}{\sqrt{1+|\nabla F^{+}({\bf x})|^{2}}}, \ \ \ {\bf x}\in\Omega(K,l), 
\]
extends continuously to a function on $\Omega(K,l)\cup K$  and $\mu(x,k(x))\cdot \nu(x)=1$  for $x\in [-c,c].$  
Now let $V$  be an open subset of $\Omega$  with $C^{2,\lambda}$  boundary such that $\{(x,-\psi(x)) : x\in [a,b]\}\subset \partial V$  and 
$\partial V\cap \left(\partial\Omega(K,l)\setminus K\right)=\emptyset$  and then let 
$U=\{(x,-y) : (x,y)\in V\}$  and $h(x,y)=F^{+}(x,-y)$  for $(x,y)\in\overline{U}.$  \qed

\begin{lem}
\label{Prop2}
Suppose $\Omega$  is a locally Lipschitz  domain in $\Real^{2},$  $\Gamma$  is a $C^{2,\lambda}$  open subset of $\partial\Omega$  for some $\lambda\in (0,1),$  
and the  curvature $\Lambda$  of $\Gamma$  is either negative or nonpositive, vanishes only at finite number of points of $\Gamma$  and (\ref{c-bound}) 
holds at each such point.  Let ${\bf y}\in\Gamma.$ 
Then there exist $\delta>0$  and upper and lower Bernstein pairs $\left(U^{\pm},\psi^{\pm}\right)$  for $(\Gamma\cap B_{\delta}({\bf y}),0).$     
\end{lem}
\vspace{3mm} 

\noindent {\bf Proof:} 
Let $\Omega\subset\Real^{2}$  be an open set,  $\Gamma\subset\partial\Omega$  be a $C^{2,\lambda}$  curve and ${\bf y}\in \Gamma$  
be a point at which we wish to have upper and lower Bernstein pairs for $H\equiv 0.$  
Notice that (\ref{c-bound}) implies that the ``curvature'' condition in Lemma \ref{Prop1} (i.e. after rotating $\Omega,$  the condition that if 
$\Lambda(\bar x,\psi(\bar x))=0,$  then $\psi''(x)\le -C_{1}|x-\bar x|^{\lambda}$  when $x$  is near $\bar x$)  is satisfied.  
Choose a neighborood $V$  of ${\bf y}$  and a rigid motion $\zeta:\Real^{2}\to\Real^{2}$  such that $\Sigma\myeq V\cap \partial\Omega \subset\Gamma$ 
and the curve $\zeta\left(\Sigma\right)$  satisfies the hypotheses of Lemma \ref{Prop1}.
Let $U$  and $h$  be as given in the conclusion of Lemma \ref{Prop1}.   
Then $\left(\zeta^{-1}(U),h\circ\zeta\right)$  will be an upper Bernstein pair for $\Sigma$  and $H\equiv 0$  and 
$\left(\zeta^{-1}(U),-h\circ\zeta\right)$  will be a lower Bernstein pair for $\Sigma$  and $H\equiv 0.$  \qed 
\vspace{3mm} 

\begin{lem}
\label{Prop3}
Suppose $\Omega$  is a $C^{2,\lambda}$  domain in $\Real^{2}$  for some $\lambda\in (0,1),$  ${\bf y}\in\partial\Omega$  and 
$\Lambda({\bf y})< -2|H({\bf y})|,$  where $\Lambda({\bf y})$  denotes the curvature of $\partial\Omega$  at ${\bf y}.$    
Then there exist $\delta>0$  and upper and lower Bernstein pairs $\left(U^{\pm},\psi^{\pm}\right)$  for $(\Gamma,H),$   where 
$\Gamma=B_{\delta}({\bf y})\cap {\partial}\Omega.$  
\end{lem}
\vspace{3mm} 

\noindent {\bf Proof:} We may assume $H\not\equiv 0$  in any neighborhood of ${\bf y}$  since Lemma \ref{Prop2} covers this case.  
There exists $\delta_{1}>0$  such that  $\Lambda({\bf x})< -2|H({\bf x})|$  for each ${\bf x}\in\partial\Omega \cap B_{\delta_{1}}({\bf y}).$   
There exists a $\delta_{2}\in (0,\delta_{1}/2)$  such that 
\[
\Lambda_{0}\myeq \sup\{\Lambda({\bf x}) : {\bf x}\in\partial\Omega\cap B_{\delta_{2}}({\bf y})\} < 
\inf\{-2|H({\bf x})| : {\bf x}\in\overline{\Omega}\cap B_{\delta_{2}}({\bf y})\} \myeq -2H_{0}.
\]
Let $z=h(\hat r),$    be a unduloid surface (see for example \cite{FinnMoon,LS1}) defined on the annulus 
${\mathcal A}= {\mathcal A}({\bf q})\myeq \{ {\bf x}\in \Real^{2} : r_{1}\le \hat r({\bf x})\le r_{2} \}$   with constant mean curvature $-H_{0}<0$  which 
becomes vertical at $\hat r=r_{1},r_{2}$   (with $Th\cdot\nu=1$  on $\hat r({\bf x})=r_{1}$), where $\hat r=|{\bf x}-{\bf p}|,$  
${\bf p}={\bf p}({\bf q})={\bf q}+r_{1}\nu({\bf q})$  for ${\bf q}\in \partial\Omega\cap B_{\delta_{2}}({\bf y}),$    
$r_{1}=\frac{1-\sqrt{1+4c_{0}H_{0}}}{2H_{0}},$   $r_{2}=\frac{1+\sqrt{1+4c_{0}H_{0}}}{2H_{0}},$  and $c_{0}\in\left(-\frac{1}{4H_{0}},0\right)$  is arbitrary.   
Now ${\mathcal A}({\bf q})$  touches $\partial\Omega$  at ${\bf q}$  and there exists $0<\delta_{3}<\delta_{2}/2$  such that 
$\{ {\bf x}\in B_{\delta_{3}}({\bf y}) : |{\bf x}-{\bf p}({\bf q})|=r_{1}\} \subset \Omega\cup \{ {\bf q} \}$  and 
$h'(r)\le 0$  if $r_{1}<r<r_{1}+\delta_{3}.$  

Set $z_{1}=h(r_{1})-h(r_{2})$  and $z_{2}>z_{1}.$  
Let $\psi^{+}$  be the variational solution of the Dirichlet problem $N\psi^{+}=-2H_{0}$  in $U=\Omega\cap B_{\delta_{3}}({\bf y}),$  
$\psi^{+}=z_{2}$  on $\Gamma=B_{\delta_{3}}({\bf y})\cap\partial\Omega$  and $\psi^{+}=0$  on $\partial U\setminus\Gamma.$
By \cite[Theorem 5.1]{FinnBook}, we see that $\psi^{+}\le z_{1}<z_{2}$  on $\Gamma$  and so $T\psi^{+}\cdot\nu=1$  on $\Gamma;$  
this follows, for example, from \cite{Bour,Lin}.  
Now define $\psi^{-}=-\psi^{+}$  in $U.$  \qed 
\vspace{3mm}

\begin{rem}
\label{Remark1} 
Suppose $\Omega$  is a $C^{2,\lambda}$  domain in $\Real^{2}$  for some $\lambda\in (0,1),$  ${\bf y}\in\partial\Omega$  and 
$\Lambda({\bf y})< 2|H({\bf y})|,$  where $\Lambda({\bf y})$  denotes the curvature of $\partial\Omega$  at ${\bf y}.$    
If $H$  is non-negative  in $U\cap\Omega$  for some neighborhood $U$  of ${\bf y},$  then the argument which establishes 
\cite[Corollary 14.13]{GT}  and boundary regularity results (e.g. \cite{Bour,Lin}) imply that there exist $\delta>0$  and an upper 
Bernstein pair $\left(U^{+},\psi^{+}\right)$  for $(\Gamma,H),$  where $\Gamma=B_{\delta}({\bf y})\cap {\partial}\Omega.$  
If $H$  is non-positive  in $U\cap\Omega$  for some neighborhood $U$  of ${\bf y},$  then there exist $\delta>0$  and a lower  
Bernstein pair $\left(U^{-},\psi^{-}\right)$  for $(\Gamma,H),$  where $\Gamma=B_{\delta}({\bf y})\cap {\partial}\Omega.$   
\end{rem}
\vspace{3mm}

\noindent {\bf Proof of Theorem \ref{Thm3}:}  
We note, as in \cite{LS1}, that the conclusion of Theorem \ref{Thm3} is a local one and so, for small $\delta>0,$  we can replace $\Omega$  by a $C^{2,\lambda}$  
set $\Omega^{*}$  such that $\Omega\cap B_{\delta}({\bf y}) = \Omega^{*}\cap B_{\delta}({\bf y})$  and $\Omega^{*}\subset  B_{2\delta}({\bf y}).$ 
We may assume $\Omega$  is a bounded domain.  Set  $S_{0} = \{ ({\bf x},f({\bf x})) : {\bf x} \in \Omega \}.$
From the calculation on page 170 of \cite{LS1},  we see that the area of $S_{0}$  is finite; let $M_{0}$  denote this area. 
For $\delta\in (0,1),$  set 
\[
p(\delta) = \sqrt{\frac{8\pi M_{0}}{\ln\left(\frac{1}{\delta}\right)}}.
\]
Let $E= \{ (u,v) : u^{2}+v^{2}<1 \}.$ 
As in \cite{EL1986A,LS1}, there is a parametric description of the surface $S_{0},$  
\begin{equation}
\label{PARAMETRIC}
Y(u,v) = (a(u,v),b(u,v),c(u,v)) \in C^{2}(E:{\Real}^{3}), 
\end{equation}
which has the following properties:
 
\noindent $\left(a_{1}\right)$  $Y$ is a diffeomorphism of $E$ onto $S_{0}$.
 
\noindent $\left(a_{2}\right)$
Set $G(u,v)=(a(u,v),b(u,v)),$  $(u,v)\in E.$   Then $G \in C^{0}(\overline{E} : {\Real}^{2}).$  

\noindent $\left(a_{3}\right)$  
Set $\sigma({\bf y})=G^{-1}\left(\partial \Omega\setminus \{ {\bf y} \}\right);$  
then $\sigma({\bf y})$ is a connected arc of $\partial E$  and $Y$ maps $\sigma({\bf y})$  onto $\partial \Omega\setminus \{ {\bf y} \}.$  
We may assume the endpoints of $\sigma({\bf y})$  are ${\bf o}_{1}({\bf y})$  and ${\bf o}_{2}({\bf y}).$  
(Note that ${\bf o}_{1}({\bf y})$  and ${\bf o}_{2}({\bf y})$  are not assumed to be distinct.)
 
\noindent $\left(a_{4}\right)$
$Y$ is conformal on $E$: $Y_{u} \cdot Y_{v} = 0, Y_{u}\cdot Y_{u} = Y_{v}\cdot Y_{v}$
on $E$.
 
\noindent $\left(a_{5}\right)$
$\triangle Y := Y_{uu} + Y_{vv} = H\left(Y\right)Y_{u} \times Y_{v}$  on $E$.
\vspace{2mm} 

\noindent From Lemma \ref{Prop2}  when $H\equiv 0$  and Lemma \ref{Prop3}  when $H$  satisfies (\ref{(16.60)}), we see that   
upper and lower Bernstein pairs $\left(U^{\pm},\psi^{\pm}\right)$  for $(\Gamma,H)$  exist.  
Notice that for each $C\in\Real,$   $Q(\psi^{+}+C)=Q(\psi^{+})\le 0$  on $\Omega\cap U^{+}$ 
and $Q(\psi^{-}+C)=Q(\psi^{-})\ge 0$  on $\Omega\cap U^{-},$  and so 
\begin{equation}
\label{Barrier+}
N(\psi^{+}+C)({\bf x}) \le 2H({\bf x},f({\bf x}))=Nf({\bf x}) \ \ {\rm for} \ \ {\bf x}\in\Omega\cap U^{+}, 
\end{equation}
and 
\begin{equation}
\label{Barrier-}
N(\psi^{-}+C)({\bf x}) \ge 2H({\bf x},f({\bf x}))=Nf({\bf x}) \ \ {\rm for} \ \ {\bf x}\in\Omega\cap U^{-}.
\end{equation}
Let $q$  denote a modulus of continuity for $\psi^{+}$  and $\psi^{-}.$  

Let $\zeta({\bf y})=\partial E\setminus\sigma({\bf y});$  then $G(\zeta({\bf y}))=\{{\bf y}\}$  and ${\bf o}_{1}({\bf y})$  and ${\bf o}_{2}({\bf y})$  
are the endpoints of $\zeta({\bf y}).$  
There exists a $\delta_{1}>0$  such that if ${\bf w}\in E$  and ${\rm dist}\left({\bf w}, \zeta({\bf y})\right)\le 2\delta_{1},$  
then $G({\bf w})\in U^{+}\cap U^{-}.$  
Now $T\psi^{\pm}\cdot \nu=\pm 1$  (in the sense of \cite{CF:74a}) almost everywhere on an open subset $\Upsilon^{\pm}$  of $\partial U^{\pm}$  
which contains $\overline{\Gamma};$  there exists a $\delta_{2}>0$  such that 
$\left(\partial U^{\pm} \setminus \Upsilon^{\pm}\right) \cap \{{\bf x}\in \Real^{2} : |{\bf x}-{\bf y}|\le 2p(\delta_{2})\}=\emptyset.$  
Set $\delta^{*}=\min\{\delta_{1},\delta_{2}\}$  and 
\[
V^{*}= \{ {\bf w}\in E : {\rm dist}({\bf w},\zeta({\bf y}))<\delta^{*} \}.
\]
Notice if ${\bf w}\in V^{*},$  then $G({\bf w})\in U^{+}\cap U^{-}.$
\vspace{3mm} 

\noindent {\bf Claim:}   $Y$  is uniformly continuous on $V^{*}$  and so extends to a continuous function on $\overline{V^{*}}.$   
\vspace{3mm} 

\noindent {\bf Pf:}  Let $\epsilon>0.$  Choose $\delta\in \left(0,\left(\delta^{*}\right)^{2}\right)$  such that  $p(\delta)+2q(p(\delta))<\epsilon.$  
Let ${\bf w}_{1},{\bf w}_{2}\in V^{*}$  with $|{\bf w}_{1}-{\bf w}_{2}|<\delta;$  then 
$G({\bf w}_{1}), G({\bf w}_{2})\in U^{+}_{1} \cap U^{-}.$
Set $C_{r}({\bf w}) = \{ {\bf u} \in E : |{\bf u} - {\bf w}| = r \}$  and  $B_{r}({\bf w}) = \{ {\bf u} \in E : |{\bf u} - {\bf w}| < r \}.$
From the Courant-Lebesgue Lemma (e.g. Lemma $3.1$ in \cite{Cour:50}), we see that there exists $\rho=\rho(\delta)\in \left(\delta,\sqrt{\delta}\right)$  
such that the arclength $l_{\rho}({\bf w}_{1})$  of $Y(C_{\rho}({\bf w}_{1}))$  is less than $p(\delta).$  
Notice that ${\bf w}_{2}\in B_{\rho(\delta)}({\bf w}_{1}).$   Let 
\[
k(\delta)({\bf w}_{1})= \inf_{{\bf u}\in C_{\rho(\delta)}({\bf w}_{1})}c({\bf u}) = \inf_{ {\bf x}\in G(C_{\rho(\delta)}({\bf w}_{1})) } f({\bf x})
\]
and 
\[
m(\delta)({\bf w}_{1})= \sup_{{\bf u}\in C_{\rho(\delta)}({\bf w}_{1})}c({\bf u}) = \sup_{ {\bf x}\in G(C_{\rho(\delta)}({\bf w}_{1})) } f({\bf x});
\]
then  
\[
m(\delta)({\bf w}_{1})-k(\delta)({\bf w}_{1})\le l_{\rho} < p(\delta).
\]
Fix ${\bf x}_{0}\in C^{\prime}_{\rho(\delta)}({\bf w}_{1}).$   
Set 
\[
C^{+}=\inf_{{\bf x}\in U^{+}\cap C^{\prime}_{\rho(\delta)}({\bf w}_{1}) } \psi^{+}({\bf x}) \ \ {\rm and} \ \ 
C^{-}=\sup_{{\bf x}\in U^{-}\cap C^{\prime}_{\rho(\delta)}({\bf w}_{1}) } \psi^{-}({\bf x}).
\]
Then $\psi^{+}-C^{+}\ge 0$  on $U^{+}\cap C^{\prime}_{\rho(\delta)}({\bf w}_{1})$  and 
$\psi^{-}-C^{-}\le 0$  on $U^{-}\cap C^{\prime}_{\rho(\delta)}({\bf w}_{1}).$ 
Therefore, for ${\bf x}\in U^{+}\cap U^{-}\cap C^{\prime}_{\rho(\delta)}({\bf w}_{1}),$  we have 
\[
k(\delta)({\bf w}_{1})+\left(\psi^{-}({\bf x})-C^{-}\right) \le f({\bf x}) \le m(\delta)({\bf w}_{1})+\left(\psi^{+}({\bf x})-C^{+}\right). 
\] 
Set 
\[
b^{+}({\bf x})= m(\delta)({\bf w}_{1})+\left(\psi^{+}({\bf x})-C^{+}\right) \ \ \ \ 
{\rm for} \ \ {\bf x}\in U^{+}\cap \overline{G\left(B_{\rho(\delta)}({\bf w}_{1})\right)}
\]
and 
\[
b^{-}({\bf x})= k(\delta)({\bf w}_{1})+\left(\psi^{-}({\bf x})-C^{-}\right) \ \ \ \ 
{\rm for} \ \ {\bf x}\in  U^{-} \cap \overline{G\left(B_{\rho(\delta)}({\bf w}_{1})\right)}.
\]
Now $\rho(\delta)<\sqrt{\delta}<\delta^{*}\le\delta_{2};$  notice that if ${\bf w}\in \overline{B_{\rho(\delta)}({\bf w}_{1})},$  then 
$|{\bf w}-{\bf w}_{1}|<\delta_{2}$  and  $|G({\bf w})-{\bf y}|<2p(\delta_{2})$  and thus if 
${\bf x}\in \overline{G\left(B_{\rho(\delta)}({\bf w}_{1})\right)} \cap \partial U^{\pm},$  then ${\bf x}\in \Upsilon^{\pm}.$  
From (\ref{Barrier+}) and (\ref{Barrier-}), the facts that $b^{-}\le f$  on $U^{-} \cap C^{\prime}_{\rho(\delta)}({\bf w}_{1})$  
and $f\le b^{+}$  on $U^{+} \cap C^{\prime}_{\rho(\delta)}({\bf w}_{1})$  and the general comparison principle (\cite[Theorem 5.1]{FinnBook}), 
we have  
\begin{equation}
\label{uniformA}
b^{-}\le f\ {\rm on} \ U^{-} \cap \overline{G\left(B_{\rho(\delta)}({\bf w}_{1})\right)} 
\end{equation}
and 
\begin{equation}
\label{uniformB}
f\le b^{+}\ {\rm on} \ U^{+}\cap \overline{G\left(B_{\rho(\delta)}({\bf w}_{1})\right)}.
\end{equation}

\noindent Since the diameter of $G\left(B_{\rho(\delta)}({\bf w}_{1})\right)\le p(\delta),$  we have 
$\left|\psi^{\pm}({\bf x})-C^{\pm}\right|\le q(p(\delta))$  for ${\bf x}\in U^{\pm}\cap G\left(B_{\rho(\delta)}({\bf w}_{1})\right).$  
Thus, whenever ${\bf x}_{1},{\bf x}_{2}\in \overline{G\left(B_{\rho(\delta)}({\bf w}_{1})\right)},$  we have 
${\bf x}_{1},{\bf x}_{2} \in U^{+}\cap U^{-}.$  
Since $c({\bf w})=f\left(G({\bf w})\right),$  $G({\bf w}_{1})\in U^{+}\cap U^{-}$  and $G({\bf w}_{2})\in U^{+}\cap U^{-},$   we have 
\[
b^{-}\left(G({\bf w}_{1})\right)-b^{+}\left(G({\bf w}_{2})\right) \le c({\bf w}_{1})-c({\bf w}_{2}) \le 
b^{+}\left(G({\bf w}_{1})\right)-b^{-}\left(G({\bf w}_{2})\right)
\]
or 
\[
-\left[m(\delta)({\bf w}_{1})-k(\delta)({\bf w}_{1}) +\left(\psi^{+}(G({\bf w}_{2}))-C^{+}\right) - \left(\psi^{-}(G({\bf w}_{1}))+C^{-}\right)\right]
\]
\[
\le c({\bf w}_{1})-c({\bf w}_{2}) \le 
\]
\[
\left[m(\delta)({\bf w}_{1})-k(\delta)({\bf w}_{1}) +\left(\psi^{+}(G({\bf w}_{1}))-C^{+}\right) - \left(\psi^{-}(G({\bf w}_{2}))+C^{-}\right)\right].
\]
Since $|\psi^{\pm}(G({\bf w}))-C^{\pm}| \le q(p(\delta))$  for ${\bf w}\in B_{\rho(\delta)}({\bf w}_{1})\cap U^{\pm},$  we have   
\[
|c({\bf w}_{1})-c({\bf w}_{2})| \le  p(\delta)+2q(p(\delta))<\epsilon.
\]
Thus $c$  is uniformly continuous on  $V^{*}$  and, since $G \in C^{0}(\overline{E} : {\Real}^{2}),$  we see that  $Y$  is uniformly continuous on $V^{*}.$  
Therefore $Y$  extends to a continuous function, still denote $Y,$  on $\overline{V^{*}}.$ \qed
\vspace{3mm} 

Notice that 
\[
\lim_{\Gamma_{1}\ni {\bf x}\to {\bf y}} f({\bf x})=\lim_{\partial E\ni {\bf w}\to {\bf o}_{1}({\bf y})} c({\bf w})=c({\bf o}_{1}({\bf y}))
\]
and 
\[
\lim_{\Gamma_{2}\ni {\bf x}\to {\bf y}} f({\bf x})=\lim_{\partial E\ni {\bf w}\to {\bf o}_{2}({\bf y})} c({\bf w})=c({\bf o}_{2}({\bf y}))
\]
and so, with $z_{1}=c({\bf o}_{1}({\bf y}))$  and $z_{2}=c({\bf o}_{2}({\bf y})),$  we see that (\ref{Apple}) holds.
  
Then the limits 
\begin{equation}
\label{Apple}
\lim_{\Gamma_{1}\ni {\bf x}\to {\bf y}} f({\bf x})=z_{1} \ {\rm and} \ 
\lim_{\Gamma_{2}\ni {\bf x}\to {\bf y}} f({\bf x})=z_{2}
\end{equation}
exist,

Now we  need to consider two cases: 
\begin{itemize}
\item[$\left(A\right)$] ${\bf o}_{1}({\bf y})= {\bf o}_{2}({\bf y}).$
\item[$\left(B\right)$] ${\bf o}_{1}({\bf y})\neq {\bf o}_{2}({\bf y}).$
\end{itemize}
These correspond to Cases 5 and 3 respectively in Step 1 of the proof of \cite[Theorem 1]{LS1}. 
\vspace{2mm} 

\noindent {\bf Case (A):}  Suppose ${\bf o}_{1}({\bf y})= {\bf o}_{2}({\bf y});$  set ${\bf o}={\bf o}_{1}({\bf y})= {\bf o}_{2}({\bf y}).$  
Then $f$  extends to a function in $C^{0}\left(\Omega \cup \{ {\bf y}\}\right)$  and case (i) of Theorem \ref{Thm3} holds. 
\vspace{2mm} 

\noindent  {\bf Pf:} Notice that $G$  is a bijection of $E\cup \{{\bf o}\}$  and $\Omega\cup \{{\bf y}\}.$  Thus we may define $f=c\circ G^{-1},$  
so $f\left(G({\bf w})\right)=c({\bf w})$  for ${\bf w}\in E\cup \{{\bf o}\};$  this extends $f$  to a function defined on $\Omega\cup \{{\bf y}\}.$  
Let $\{\delta_{i}\}$  be a decreasing sequence of positive numbers converging to zero and consider the sequence of open sets 
$\{G(B_{\rho(i)}({\bf o}))\}$  in $\Omega,$  where $\rho(i)=\rho(\delta_{i}({\bf o})).$  
Now ${\bf y}\notin G(C_{\rho(i)}({\bf o}))$  and so there exist $\sigma_{i}>0$  such that 
\[
P(i)  =\{ {\bf x}\in\Omega : |{\bf x}-{\bf y}|<\sigma_{i}\} \subset G(B_{\rho(i)}({\bf o}))
\]
for each $i\in\Natural.$  Thus if  ${\bf x}\in P(i),$  we have $|f({\bf x}) - f({\bf y})|<p(\delta_{i})+2q(p(\delta_{i})).$  
The continuity of $f$  at ${\bf y}$  follows from this.  \qed 
\vspace{3mm} 

\noindent {\bf Case (B):}  Suppose ${\bf o}_{1}({\bf y})\neq {\bf o}_{2}({\bf y}).$   Then case (ii) of Theorem \ref{Thm3} holds.
\vspace{2mm} 

\noindent  {\bf Pf:}
As at the end of Step 1 of the proof of \cite[Theorem 1]{LS1}, we define $X:B\to\Real^{3}$  by $X=Y\circ g$ and $K:B\to\Real^{2}$  
by $K=G\circ g,$  where $B=\{(u,v)\in\Real^{2} : u^{2}+v^{2}<1, \ v>0\}$  and $g:\overline{B}\to \overline{E}$  is either a conformal or an indirectly 
conformal (or anticonformal) map from $\overline{B}$  onto $\overline{E}$  such that  $g(1,0)= {\bf o}_{1}({\bf y}),$   
$g(-1,0)= {\bf o}_{2}({\bf y})$  and $g(u,0)\in  {\bf o}_{1}({\bf y}){\bf o}_{2}({\bf y})$  for each $u\in [-1,1],$  where 
${\bf ab}$  denotes the (appropriate) choice of arc in $\partial E$  with ${\bf a}$  and ${\bf b}$  as endpoints.  

Notice that $K(u,0)={\bf y}$  for $u\in [-1,1].$  Set $x=a\circ g,$  $y=b\circ g$  and $z=c\circ g,$  so that $X(u,v)=(x(u,v),y(u,v),z(u,v))$  
for $(u,v)\in B.$  Now, from Step 2 of the proof of \cite[Theorem 1]{LS1}, 
\[
X\in C^{0}\left(\overline{B}:\Real^{3}\right)\cap C^{1,\iota}\left(B\cup\{(u,0):-1<u<1\}:\Real^{3}\right)
\]
for some $\iota\in (0,1)$  and $X(u,0)=({\bf y},z(u,0))$  cannot be constant on any nondegenerate interval in $[-1,1].$  
Define $\Theta(u)= {\rm arg}\left( x_{v}(u,0)+iy_{v}(u,0) \right).$  From equation \cite[(12)]{LS1}, we see that 
\[
\alpha_{1}=\lim_{u\downarrow -1} \Theta(u) \ \ \ \ {\rm and} \ \ \ \  \alpha_{2}=\lim_{u\uparrow 1} \Theta(u);
\]
here $\alpha_{1}<\alpha_{2}.$  
As in Steps 2-5 of the proof of \cite[Theorem 1]{LS1}, we see that $Rf(\theta)$  exists when $\theta\in \left(\alpha_{1},\alpha_{2}\right),$ 
\[
\overline{G^{-1}\left( L(\alpha_{2}) \right)} \cap \partial E = \{ {\bf o}_{1}({\bf y}) \} \    (\& \ 
\overline{K^{-1}\left( L(\alpha_{2}) \right)} \cap \partial B = \{ (1,0) \})   \  {\rm when}  \  \alpha_{2}<\beta({\bf y})
\]
\[
\overline{G^{-1}\left( L(\alpha_{1}) \right)} \cap \partial E = \{ {\bf o}_{2}({\bf y}) \} \ (\& \ 
\overline{K^{-1}\left( L(\alpha_{1}) \right)} \cap \partial B = \{ (-1,0) \})  \  {\rm when} \ \alpha_{1}>\alpha({\bf y})
\]
where  $L(\theta)= \{{\bf y}+(r\cos(\theta),r\sin(\theta))\in \Omega : 0<r<\delta^{*} \},$
and we see that  $Rf$ is strictly increasing or strictly decreasing on $(\alpha_{1}, \alpha_{2}).$    
We may argue as in Case A to see that $f$  is uniformly continuous on 
\[
\Omega^{+} =\{{\bf y}+(r\cos(\theta),r\sin(\theta))\in \Omega : 0<r<\delta, \alpha_{2}\le \theta<\beta({\bf y})+\epsilon\}
\]
and $f$  is uniformly continuous on 
\[
\Omega^{-} =\{{\bf y}+(r\cos(\theta),r\sin(\theta))\in \Omega : 0<r<\delta, \alpha({\bf y})-\epsilon< \theta\le\alpha_{1}\}
\]
for some small $\epsilon>0$  and $\delta>0,$  since  $G$  is a bijection of $E\cup \{{\bf o}_{1}({\bf y})\}$  and $\Omega\cup \{{\bf y}\}$  
and a bijection of $E\cup \{{\bf o}_{2}({\bf y})\}$  and $\Omega\cup \{{\bf y}\}.$  
Theorem \ref{Thm3} then follows, as in \cite{NoraKirk1}, from Steps 2-5 of the proof of \cite[Theorem 1]{LS1}.  
\qed
\vspace{3mm}

\noindent {\bf Proof of Theorem \ref{Thm2}:}   
From Theorem \ref{Thm3}, we see that the radial limits $Rf(\theta,{\bf y})$  exist for each $\theta\in [\alpha({\bf y}),\beta({\bf y})].$  
Set $z_{1}=Rf(\alpha({\bf y}),{\bf y}),$  $z_{2}=Rf(\beta({\bf y}),{\bf y})$  and $z_{3}=\phi({\bf y}).$  
If $z_{1}=z_{2},$  then case (i) of Theorem \ref{Thm3} holds.  
(If $f$  is symmetric with respect to a line through ${\bf y},$  then $z_{1}=z_{2}$  and we are done.)

Suppose otherwise that $z_{1}\neq z_{2};$  we may assume that $z_{1}<z_{3}$  and $z_{1}<z_{2}.$    
Then there exist $\alpha_{1}, \alpha_{2}\in [\alpha({\bf y}),\beta({\bf y})]$  with $\alpha_{1}<\alpha_{2}$  such that 
\[
Rf(\theta,{\bf y}) \ \ {\rm is} 
\left\{ \begin{array}{ccc} 
{\rm constant}(=z_{1}) & {\rm for} & \alpha({\bf y})\le\theta\le\alpha_{1}\\ 
{\rm strictly \ increasing} & {\rm for} & \alpha_{1}\le\theta\le\alpha_{2}\\
{\rm constant}(=z_{2}) & {\rm for} & \alpha_{2}\le\theta\le\beta({\bf y}). \\ 
\end{array}
\right.
\]
From Theorem \ref{Thm3}, we see that $Rf(\theta,{\bf y})$  exists for each $y\in \Gamma$  and $\theta\in [\alpha({\bf y}),\beta({\bf y})]$   
and $f$  is continuous on $\Omega\cup\Gamma\setminus \Upsilon$  for some countable subset $\Upsilon$  of $\Gamma.$  
Let $z_{0}\in \left(z_{1},\min\{z_{2},z_{3}\}\right)$  and $\theta_{0}\in (\alpha_{1},\alpha_{2})$  satisfy $Rf(\theta_{0},{\bf y})=z_{0}.$  
Let $C_{0}\subset\Omega$  be the $z_{0}-$level curve of $f$  which has ${\bf y}$  and a point ${\bf y}_{0}\in \partial\Omega\setminus \{{\bf y}\}$  as endpoints. 
Let ${\bf y}_{1}\in \Gamma_{1}({\bf y})\cap\Gamma\setminus \Upsilon$  and ${\bf y}_{2}\in C_{0}$  such that the (open) line segment $L$  joining
${\bf y}_{1}$  and ${\bf y}_{2}$  is entirely contained in $\Omega.$  
Let $M<\min\{z_{1},\inf_{L} f\}.$  
Let $\Pi$  be the plane containing $({\bf y},z_{0})$  and $L\times \{ M\}$  and let $h$  be the affine function on $\Real^{2}$  whose graph is $\Pi.$   
Let $\Omega_{0}$  be the component of $\Omega\setminus \left(C_{0}\cup L\right)$  whose closure contains $B_{\delta}({\bf y})\cap \Gamma_{1}({\bf y})$  
for small $\delta>0.$  
By first choosing ${\bf y}_{2}$  sufficiently near ${\bf y}$  and then choosing ${\bf y}_{1}$  sufficiently near ${\bf y},$  we may assume that 
${\bf y}$  is the furthest point on $C_{0}$  between ${\bf y}$  and ${\bf y}_{2}$  away from $L.$  
Then $h=M<f$  on $L,$  $h\le z_{0}=f$  on the portion of $C_{0}$  between  ${\bf y}$  and ${\bf y}_{2},$ 
and $h({\bf y})=z_{0}>Rf(\theta,{\bf y})$  for each $\theta \in [\alpha({\bf y}),\theta_{0}).$  
Thus $h\le f$  on $\Omega\cap\partial\Omega_{0}$  and $h>f$  on $B_{\delta}({\bf y})\cap \partial\Omega_{0}$  for some sufficiently small $\delta>0.$

\begin{figure}[ht]
\centering
\includegraphics{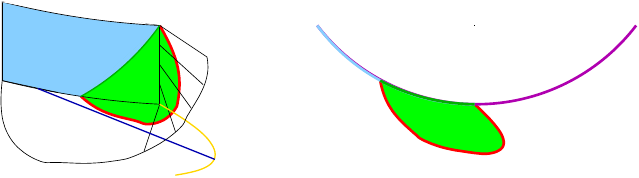}
\caption{Side View of $\Pi\cap \left(\Omega\times\Real\right)$ (left) and $\Omega_{1}$ (right)  \label{FOUR_F}}
\end{figure}

\begin{figure}[ht]
\centering
\includegraphics{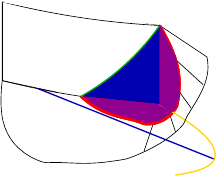}
\caption{Graph of $f$  over $\Omega_{1}$  (purple) and part of $\partial\Omega\times\Real$  (blue)  \label{FOUR_FC}}
\end{figure}

Then there is a curve $C\subset \Omega_{0}$  on which $f=h$   whose endpoints are ${\bf y}_{3}$  and ${\bf y},$  for some 
${\bf y}_{3}\in \Gamma_{1}({\bf y})$  between ${\bf y}_{1}$  and ${\bf y},$  such that $h>f$  in $\Omega_{1},$  where $\Omega_{1}\subset\Omega_{0}$  
is the open set bounded by $C$  and the portion of $\Gamma_{1}({\bf y})$  between ${\bf y}$  and ${\bf y}_{3}.$    
(In Figure \ref{FOUR_F}, on the left,  $\{\left({\bf x},h\left({\bf x}\right)\right): {\bf x}\in C\}$  is in red, $L$  is in dark blue, $C_{0}$  is in yellow, and 
the light blue region is a portion of $\partial_{{\bf y}}^{1}\Omega\times\Real,$  and, on the right, $\Omega_{1}$  is in light green and 
$\Gamma_{2}({\bf y})$  is in magenta.) 
Now let $g\in C^{2}(\Omega)$  be defined by $g=f$  on $\Omega\setminus\overline{\Omega_{1}}$  and $g=h$  on $\Omega_{1}$  and observe that 
$J(g)<J(f).$  (The functional $J(f)$  includes the area of the blue surface in Figure \ref{FOUR_FC}, which is a subset of $\partial\Omega\times\Real,$ 
and the area of the purple surface of this figure, which is the subset of the graph of $f$  over $\Omega_{1}$  while $J(g)$  does not include the areas 
of the blue and purple surfaces and instead includes the area of the green surface on the left side of Figure \ref{FOUR_F}, 
which is the portion of the plane $\Pi$  over $\Omega_{1}.$)
This contradicts the fact that $f$  minimizes $J.$  Thus it must be the case that $z_{1}=z_{2},$  case (i) of Theorem \ref{Thm3} holds 
and $f$  is continuous at ${\bf y}.$  
(Notice that the set $\Omega_{1}=\{ {\bf x}\in\Omega_{0} : h({\bf x})>f({\bf x})\}$  could be more complex but the proof is unchanged; if 
$V$  is an open set in $\Omega_{0}$  with $h<f$  in $V,$  $h=f$  on $\partial V\cap\Omega$  and $h\le f$  on $\partial V\cap \partial\Omega$  
is an ``inclusion'' in $\Omega_{1},$  it does not matter; we still set $g=h$  on $\Omega_{1}$  and $g=f$  on $\Omega\setminus\Omega_{1}.$) 
\qed  

\vspace{3mm}

\noindent {\bf Proof of Theorem \ref{Thm1}:} 
The proof of Theorem \ref{Thm1} is essentially the same as that of Theorem \ref{Thm2};  
we replace $f$  with a function $g$  and obtain a contradiction by showing $J(g)<J(f),$  where 
\[
J(u)=\int_{\Omega} \sqrt{1+|Du|^{2}} d{\bf x} + \int_{\Omega} \left( \int_{c}^{u({\bf x})} 2H({\bf x},t)  \ dt \right) d{\bf x}
+ \int_{\partial\Omega} |u-\phi| ds
\]
for $u\in BV(\Omega).$  

\begin{figure}[ht]
\centering
\includegraphics{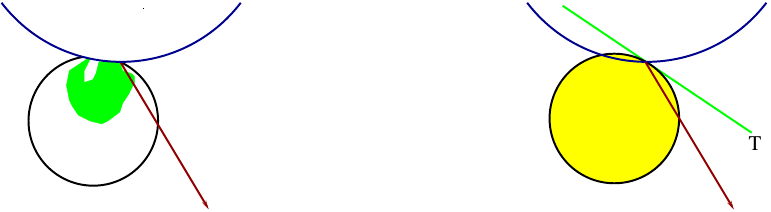}
\caption{$U$  (green region; left); $C(R)$  (black circle; right) \label{FIVE}}
\end{figure}

As in the proof of Theorem \ref{Thm2}, we assume $z_{1}<z_{3},$  $z_{1}<z_{2},$  and  $Rf(\theta,{\bf y})$  is strictly increasing on 
$(\alpha_{1},\alpha_{2})$  and constant on $[\alpha({\bf y}),\alpha_{1}]$  and $[\alpha_{2},\beta({\bf y})].$  
Let $z_{0}\in \left(z_{1},\min\{z_{2},z_{3}\}\right)$  and $\theta_{0}\in (\alpha_{1},\alpha_{2})$  satisfy $Rf(\theta_{0},{\bf y})=z_{0}.$  
Extend $H$  so $H\in C^{1}(\Real^{2}).$  
Let $R>0$  be small enough that $2R|H({\bf x})|\le 1$  for  all ${\bf x}\in B_{2R}({\bf y}).$ 
Set $\theta_{b}=(\theta_{0}+\beta({\bf y}))/2$  and $T=\{ {\bf y}+r\left(\cos\theta_{b},\sin\theta_{b}\right) : r\in\Real\},$  and let $C(R)$  
be the circle of radius $R$  which passes through ${\bf y},$  is tangent at ${\bf y}$  to the line $T,$  and intersects $\Gamma_{1}({\bf y}).$  
Let $V(R)$   be the open disk inside $C(R).$  
(In Figure \ref{FIVE}, the blue arcs are part of $\partial\Omega,$  locally $\Omega$  lies below these blue arcs, the red rays represent $\theta=\theta_{0},$  
and, on the right, $V(R)$  is the yellow disk,  $C(R)$  is the black circle, and $T$  is the green line.) 
Notice that (\ref{eq:D})-(\ref{bc:D}) in the domain $V(R)$  is solvable for all $\phi\in C^{0}(C(R)).$  

Let $h\in C^{2}\left(\overline{V(R)}\right)$  satisfy $Qh=0$  in $V(R),$  $h({\bf y})=z_{0},$  
$h < f$  on $C(R)\cap\Omega$  (recall $Rf(\theta_{b},{\bf y})>Rf(\theta_{0},{\bf y})=z_{0}$),   
and $h\le z_{0}$  on $C(R)\setminus \Omega.$  
Set $U=\{ {\bf x}\in\Omega\cap V(R) : h({\bf x})>f({\bf x}) \}$  (see Figure \ref{FIVE}) and notice that 
$\Gamma_{1}({\bf y})\cap B_{\delta}({\bf y})\subset \overline{U}\subset V(R)\cup \{{\bf y}\}$  for some $\delta>0.$      
Now define $g\in C^{0}(\Omega)$  by $g=h$  on $U$  and $g=f$  on $\Omega\setminus U.$  
Then $J(g)<J(f)$  as before in the proof of Theorem \ref{Thm2}.  \qed
\vspace{3mm}

\noindent {\bf Proof of Example \ref{Example1}:} 
By Theorem \ref{Thm1}, $f$  is continuous on $\Omega\cup \{(0,0)\}.$  Clearly $f$  is continuous at $(x,y)$  when $(x+1)^2+y^2=\cosh^{2}(1).$ 
By \cite{Simon}, $f$  is continuous at $(x,y)$  when $(x+1)^2+y^2=1$  and $(x,y)\neq(0,0).$  
The parametrization (\ref{PARAMETRIC}) of the graph of $f$  (restricted to $\Omega\setminus \{(x,0):x<0\}$) satisfies $Y\in C^{0}(\overline{E}).$  
Notice that $\zeta((0,0))=\{ {\bf o} \}$  (since $\beta((0,0))-\alpha((0,0))=\pi$  and $z_{1}=z_{2}$)  for some ${\bf o}\in \partial E.$ 
Suppose $G$  in $(a_{2})$  is not one-to-one.  Then there exists a nondegenerate arc $\zeta\subset\partial E$  such that $G(\zeta)=\{{\bf y}_{1}\}$  for some 
${\bf y}_{1}\in \partial\Omega$  and therefore $f$  is not continuous at ${\bf y}_{1},$  which is a contradiction. 
Thus $f=g\circ G^{-1}$  and so $f\in C^{0}\left(\overline{\Omega}\right).$  
(The continuity of $G^{-1}$  follows, for example, from Lemma $3.1$ in \cite{Cour:50}.) 
\qed 

\begin{rem}
\label{Remark2} 
The term ``order of non-uniformity'', used in \cite{Lau}, and the term genre, established in 1912 by Bernstein, adopted by Serrin (\cite[p. 425]{Ser:69}) 
and used in \cite{JL}, are numbers which measure the variation  from  uniform ellipticity of a quasilinear elliptic operator; 
the genre and the ``order of non-uniformity'' for the minimal surface operator are both equal to two ($2$) while  
the genre and the ``order of non-uniformity'' for uniformly elliptic operators are both equal to zero ($0$).  
We are confused by the invention for a new phrase for an existing and named phenomenon.  
\end{rem}

\end{document}